\documentclass[submission]{FPSAC2026}

\usepackage{amsfonts,mathtools,bm,bbm,array,float,mathdots,dsfont,mathrsfs,latexsym,booktabs}
\usepackage{caption,subcaption,multirow,longtable,lscape,wrapfig,tikz,tikz-cd,comment}

\usetikzlibrary{calc}
\newcommand\dyckpath[5]{
  \begin{scope}[local bounding box=#4]
    \fill[white]  (#1) rectangle +(#2,#2);
    \path[fill,red] (#1) foreach \num [count=\i from 0] in {#5}{ +(-0.5,\i+0.75) node[anchor=north]{\num} \ifnum\i>#2 circle (1pt) \fi};
    \draw[help lines] (#1) grid +(#2,#2);
    \draw[help lines] (#1)--+(#2,#2);
    \draw[line width=2pt] (#1) foreach \dir in {#3}{ -- ++(\dir*90:1)};
  \end{scope}
}

\newtheorem{theorem}{Theorem}

\newtheorem{ex}{Example}

\usepackage{lipsum}

\title[From Black Box to Bijection]{From Black Box to Bijection: Interpreting Machine Learning to Build a Zeta Map Algorithm}

\author[]{Xiaoyu Huang\thanks{\href{mailto:xiaoyu.huang@temple.edu}{xiaoyu.huang@temple.edu}.}\addressmark{1} \and Blake Jackson\thanks{\href{mailto:blake.jackson@uconn.edu}{blake.jackson@uconn.edu}}\addressmark{2} \and Kyu-Hwan Lee\thanks{\href{mailto:khlee@math.uconn.edu}{khlee@math.uconn.edu}}\addressmark{2,3}}

\address{\addressmark{1}Department of Mathematics, Temple University, Philadelphia, PA 19122, U.S.A. \\ \addressmark{2}Department of
Mathematics, University of Connecticut, Storrs, CT 06269, U.S.A. \\ \addressmark{3}Korea Institute for Advanced Study, Seoul 02455, Republic of Korea}

\received{\today}

\abstract{There is a large class of problems in algebraic combinatorics which can be distilled into the same challenge: construct an explicit combinatorial bijection. Traditionally, researchers have solved challenges like these by visually inspecting the data for patterns, formulating conjectures, and then proving them. But what is to be done if patterns fail to emerge until the data grows beyond human scale? In this paper, we propose a new workflow for discovering combinatorial bijections via machine learning. As a proof of concept, we train a transformer on paired Dyck paths and use its learned attention patterns to derive a new algorithmic description of the zeta map, which we call the \textit{Scaffolding Map}.}


\usepackage[backend=bibtex]{biblatex}
\begin{document}

\maketitle

\section{Motivation} \label{sec:intro}
The $q,t$-Catalan numbers are a family of polynomials in the variables $q$ and $t$ which refine the classical Catalan numbers \cite{haglund_q_nodate}.
The origin of these polynomials dates back to 1996 and is rooted in symmetric function theory, specifically the theory of Macdonald polynomials and the study of modules over the space of diagonal harmonics \cite{garsia_remarkable_1996}. 
More precisely, the $q,t$-Catalan number $C_n(q,t) \in \mathbb{Z}_{\geq0}[q,t]$ is the bigraded Hilbert series of the alternating component of the space of diagonal harmonics in $2n$ variables.
Here are the first three $q,t$-Catalan numbers:
\begin{center}
    \begin{tabular}{c|c}
         $n$ & $C_n(q,t)$ \\
         \midrule
         1 & 1 \\
         2 & $q+t$ \\
         3 & $q^3 + q^2t + qt + qt^2 + t^3$
    \end{tabular}
\end{center}

Since $C_n(q,t)$ are a refinement of the Catalan numbers, the search began for a combinatorial interpretation for these polynomials: a pair of \textit{statistics} on some collection of objects enumerated by the Catalan numbers that would realize the coefficients of these polynomials.
Ultimately, two formulas (discovered simultaneously) emerged, leveraging statistics on Dyck paths.
\begin{theorem}[Haglund 2003, Haiman 2000\nocite{Haiman2000_qtCatalan}]
    There exists nonnegative statistics $area, bounce,$ and $dinv$ on Dyck paths such that
    \begin{align*}
        C_n(q,t) &= \sum_{w \in Dyck(n)}  q^{area(w)}t^{bounce(w)} \\
        &= \sum_{w \in Dyck(n)} q^{dinv(w)}t^{area(w)}.
    \end{align*}
\end{theorem}
In the same paper which produced the bounce statistic, Haglund \cite{haglund_conjectured_2003} also described a bijection on Dyck paths (called the \textit{zeta map} $\zeta$) with the property that \[ (area(w), bounce(w)) = (dinv(\zeta(w)), area(\zeta(w))). \]
Moreover, there have been various combinatorial descriptions of the zeta map (and its inverse), which come in two main flavors: maps on the area vector of a Dyck path \cite{andrews_ad-nilpotent_2002,haglund_conjectured_2003} and sweep maps \cite{armstrong2015sweep,thomas2018sweeping}.

The other side of the motivation for this paper comes from machine learning (ML).  
Its power and versatility have made it a valuable instrument for mathematical exploration. 
In recent years, a growing body of work under the theme of "machine learning for mathematics" has investigated whether ML models can recognize and exploit hidden mathematical patterns. 
When trained on datasets generated from mathematical objects, ML algorithms have sometimes achieved remarkably high classification accuracies in distinguishing objects according to their invariants or other defining features. 
Representative examples include algebraic curves \cite{he2022machine_sato, he2023machine, he2025murmurations, babei2024machine} and number fields \cite{amir2023machine, he2022machine}.  
More intriguingly, interpreting what the model learns has occasionally led to new insights into the underlying mathematics. 
The recently uncovered \textit{murmuration} phenomenon in arithmetic provides a striking instance of such progress \cite{quanta, he2025murmurations}.
This line of inquiry suggests a novel paradigm for mathematical research: one begins by constructing datasets, applies machine learning techniques, analyzes the learned representations, and, from these interpretations, formulates and ultimately proves conjectures.

Combining these two compelling aspects, our main contribution is a new combinatorial algorithm for the zeta map, derived from a trained \textit{transformer} model.
\begin{theorem}[Huang, Jackson, Lee 2025+]
    The zeta map is equivalent to a scaffolding map.
\end{theorem}
The scaffolding map is reminiscent of the sweep map implementations of the zeta map; the main difference is its use of the implicit symmetry in the peaks of a Dyck path.

A broader goal of this paper is to use the zeta map as a proving ground for machine learning (ML) in combinatorics research, namely to demonstrate the ability to extract \textit{combinatorial bijections} and \textit{algorithms} from trained models.
While the precise details may differ from one application of ML in math to another, the concepts articulated in this paper are widely applicable within combinatorics research.


\section{Background on Machine Learning} \label{sec:bg}

As mentioned above, we adopt a \textit{transformer} architecture as our principal learning model and apply it to the datasets of Dyck paths obtained from the zeta map. 
The transformer model, first proposed in the landmark paper "Attention is All You Need" \cite{vaswani2017attention}, has demonstrated extraordinary effectiveness in sequence-to-sequence learning tasks, particularly in natural language processing (NLP). 
It now forms the fundamental framework for large language models (LLMs). More recently, transformers have proven successful in addressing various mathematical problems, ranging from basic computational tasks such as arithmetic \cite{charton2023learning} and eigenvalue estimation for symmetric matrices \cite{charton2021linear}, to more challenging or open problems such as constructing Lyapunov functions for dynamical stability \cite{alfarano2024global}, determining certain Euler factors of elliptic curves from the other ones \cite{babei2025learning}, and even identifying optimal solutions to longstanding combinatorial problems \cite{charton2024patternboost}.  
Given that our setting naturally involves mapping input sequences to output sequences, the transformer is an especially suitable choice of model.

\subsection{Model: Transformer Architecture}
\label{s:models}
For a comprehensive account of the transformer architecture, we refer the reader to \cite{elhage2021mathematical,vaswani2017attention}.  
Here we provide a brief overview of the essential components relevant to our implementation.

The first stage of the model involves converting numerical data into discrete symbols in a manner analogous to word tokenization in natural language processing.  
In our setting, the tokenizer constructs a vocabulary consisting of $\{0,1,bos,eos\}$ that appear in the dataset, each treated as an individual token.  
The model thus operates on symbolic sequences rather than directly on numerical values.  
Each token is subsequently mapped to a $d$-dimensional vector through a learnable embedding matrix of size $v \times d$, where $v$ denotes the vocabulary size.  
This process, achieved by multiplying the token's one-hot representation by the embedding matrix, embeds tokens into a continuous vector space that the model can manipulate effectively.

Because the transformer architecture lacks an inherent sense of sequential order, positional information must be explicitly supplied.  
To this end, positional encodings are added to the token embeddings, providing the model with knowledge of each token's location within the sequence.  
These encodings may be fixed or trainable parameters optimized jointly with the model.

The transformer can follow an encoder–decoder architecture, originally designed for sequence-to-sequence tasks such as translation. 
The encoder receives an input sequence, embeds and enriches it with positional information, and processes it through a stack of $n$ transformer layers. 
Each encoder layer consists of two main components: a multi-head self-attention mechanism and a position-wise feedforward network.
The self-attention module enables every token to attend to all others in the same sequence, thereby capturing dependencies of various lengths and types.
From each token embedding, the model derives three vectors---a query, a key, and a value---which collectively form the matrices $Q$, $K$, and $V$ for the entire sequence.
The attention mechanism compares each query with all keys to compute attention weights that quantify how strongly one token should focus on another.
This allows information to flow adaptively across the sequence, updating each token's representation according to its learned relationships with others.

Rather than relying on a single attention map, the transformer employs multiple attention heads, each designed to capture distinct types of correlations within the data.
Each head performs an independent attention computation, and the resulting outputs are concatenated and linearly projected to produce the final representation. 
Each head can thus be viewed as an independent channel through which the model exchanges and processes information. 
Following the attention modules, each sequence position is independently transformed by a one-hidden-layer feedforward network: a linear projection, a nonlinear activation, and a second linear projection.


The decoder also consists of a stack of transformer layers, but introduces an additional cross-attention mechanism between the encoder and decoder. 
Within each decoder layer, the first submodule performs masked self-attention so that each position can only attend to earlier positions in the output sequence, ensuring autoregressive generation. 
The second submodule, called cross-attention, allows the decoder to attend to the encoder's output representations.
This step integrates information from the encoded input sequence, effectively aligning the target tokens with relevant parts of the source sequence. 
From an interpretive viewpoint, cross-attention provides an explicit mapping between input and output structures.

In summary, the transformer encodes an input sequence by embedding its tokens, enriching them with positional information, and successively processing them through stacked attention and feedforward layers in the encoder. 
The decoder then generates an output sequence, guided both by its own past outputs (through self-attention) and by the encoder's contextualized representations (through cross-attention).
The final decoder outputs are projected onto the vocabulary to produce a probability distribution over possible tokens at each position, from which the most likely sequence is chosen as the model's prediction.


\section{Existing algorithms} \label{sec:exts-alg}

There are (roughly) two existing combinatorial algorithms for the zeta map: area sequence maps and sweep maps.
The sweep maps have a similar flavor to our algorithm.
The first type is as a function on the area sequence/bounce path \cite{andrews_ad-nilpotent_2002,haglund_conjectured_2003}, which can be described by the following algorithm:
\begin{enumerate}
    \item Compute the area sequence $a_1, a_2, \ldots, a_n$ of $w$ (row lengths from bottom)
    \item Set $b = \text{max}\{ a_i \}+1$ to be one more than the length of the longest row
    \item Scan the area sequence from left to right, looking for $-1$'s and $0$'s
    \begin{enumerate}
        \item For each $-1$ encountered, write a 0
        \item For each $0$ encountered, write a 1
    \end{enumerate}
    \item Repeat step 3 using the integers $k$ and $k+1$ for $k = 0, 1, \ldots, b$, where encountering the smaller number appends a 0 to the output and the larger appends a 1 to the output
\end{enumerate}

The second version of the zeta map is a sweep map composed with a reversal \cite{armstrong2015sweep,thomas2018sweeping}, which can be described by the following algorithm:
\begin{enumerate}
    \item Start with the binary sequence representation $w = (w_i)_{i=1,\ldots,2n}$ of length $2n$ for the Dyck path over the alphabet $\{0, 1\}$
    \item Reverse the word to obtain $rev(w) = (w_{n-i+1})_{i=1,\ldots,2n}$
    \item Assign levels to each position of $rev(w)$ by the following recursive process: \[ \ell_0 = 0 \qquad \ell_{i+1} = \begin{dcases}
        \ell_i + 1 & \text{if } rev(w)_i = 1 \\
        \ell_i - 1 & \text{if } rev(w)_i = 0 
    \end{dcases}\]
    \item "Sweep" the reversed word together with its level sequence from left to right repeatedly, recording the entry of the reversed word that corresponds to level $\ell = 0, -1, -2, \ldots, 4, 3, 2, 1$
\end{enumerate}
Below is an example of a sweep-style zeta map acting on a Dyck vector of semilength 8.
\begin{ex}
    Consider the Dyck path $w$ below.
    \[\begin{tikzpicture}[scale=0.5]
      \dyckpath{0,0}{8}{1,1,1,0,1,0,1,1,0,0,0,1,1,0,0,0}{dyck1}{};
    \end{tikzpicture}\]
    To implement the zeta map in the spirit of \cite{armstrong2015sweep}, there are four steps we need to take: write down the binary sequence for $w$, reverse it, calculate levels, and then sweep. To any binary sequence, a \textit{level vector} $\ell_i$ can be calculated recursively: \[ \ell_0 = 0 \qquad \ell_{i+1} = \begin{dcases}
        \ell_i + 1 & \text{if } w_i = 1 \\
        \ell_i - 1 & \text{if } w_i = 0.
    \end{dcases} \]
    Following steps 1-3, we have \begin{center}
        \begin{tabular}{*{18}{c}}
         $w$ & = & 1 & 1 & 1 & 0 & 1 & 0 & 1 & 1 & 0 & 0 & 0 & 1 & 1 & 0 & 0 & 0 \\
         reverse($w$) & = & 0 & 0 & 0 & 1 & 1 & 0 & 0 & 0 & 1 & 1 & 0 & 1 & 0 & 1 & 1 & 1 \\
         $\ell$(reverse($w$)) & = & -1 & -2 & -3 & -2 & -1 & -2 & -3 & -4 & -3 & -2 & -3 & -2 & -3 & -2 & -1 & 0 \\
    \end{tabular}
    \end{center}
    We then "sweep" the binary sequence and the corresponding level vector from left to right repeatedly, recording the entry of the sequence for levels $\ell = 0, -1, -2, -3, \ldots.$
    Doing this results in the image of the Dyck vector under the zeta map:
    \begin{center}
        \begin{tabular}{*{18}{c}}
         $\zeta(w)$ & = & 1 & 0 & 1 & 1 & 0 & 1 & 0 & 1 & 1 & 1 & 0 & 0 & 1 & 0 & 0 & 0 \\
         corresponding levels & = & 0 & -1 & -1 & -1 & -2 & -2 & -2 & -2 & -2 & -2 & -3 & -3 & -3 & -3 & -3 & -4
    \end{tabular}
    \end{center}
    On the lattice path, $\zeta(w)$ looks like a reversal followed by "sweeping" along the hyperplanes $y-x=k$ as $k$ ranges from $0$ to $-\infty$, collecting the steps crossed in level order:
    \[\begin{tikzpicture}[scale=0.5]
      \dyckpath{0,0}{8}{1,1,1,0,1,0,1,1,0,0,0,1,1,0,0,0}{dyck1}{};
      \node at (9.5,4) {$reverse$};
      \node at (9.5,3.5) {$\longrightarrow$};
      \dyckpath{11,0}{8}{0,0,0,1,1,0,0,0,1,1,0,1,0,1,1,1}{dyck2}{};
      \node at (20.5,4) {$sweep$};
      \node at (20.5,3.5) {$\longrightarrow$};
      \dyckpath{22,0}{8}{1,0,1,1,0,1,0,1,1,1,0,0,1,0,0,0}{dyck3}{};
    \end{tikzpicture}\]
\end{ex}

\section{Results} \label{sec:results}

\subsection{Experiment results}
We began by experimenting with transformer models of varying complexities. 
Using a fixed architecture with an embedding dimension of 256, four encoder layers, and four decoder layers (each with eight attention heads), the model achieved 99\% accuracy in learning the zeta map, that is, predicting the image of a Dyck path under the zeta map, for semi-lengths $n = 11, 12, 13, 14, 15, 16$\footnote{The dataset sizes are $58786, 208012, 742900, 2674440, 9694845, 35357670$, corresponding to Catalan numbers $C(n)$.}.


Having established that a transformer architecture can successfully learn the zeta map, we next reduced the model's complexity and worked with one of the simplest transformer configurations: an embedding dimension of 128, a single encoder layer, a single decoder layer, and one attention head. 
The dataset consisted of all Dyck paths of semi-length $n = 13$, totaling $742{,}900$ samples, and the resulting model contained $339{,}716$ parameters. 
The hypothesis is that there is an underlying algorithm that the transformer model is implementing, and that one can uncover it when the model is sufficiently simple. 

Again, this minimal transformer, which we name \textit{Minimal Dyck Transformer}, consistently learned the zeta map with 100\% accuracy across multiple random configurations.

\subsection{Interpretation of \textit{Minimal Dyck Transformer}}
\begin{figure}
    \centering
    \includegraphics[width=0.9\linewidth]{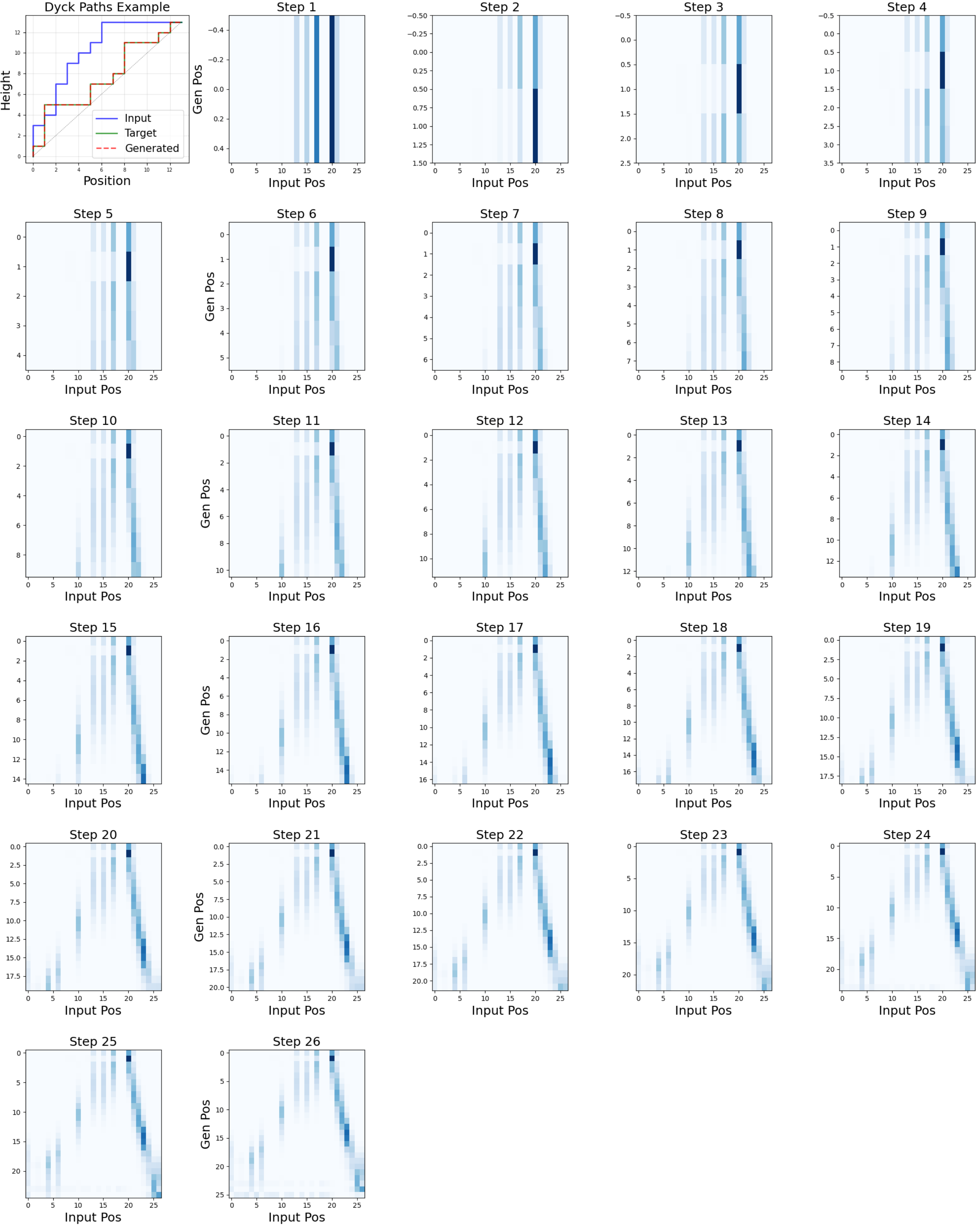}
    \caption{The visualizations show the cross-attention matrices across all generation steps for a single example, where the zeta map input is $1,1,1,0,1,0,1,1,1,0,1,1,0,1,0,1,0,1,1,0,0,0,0,0,0,0$. In the visualization for Step $k$, the entry at position $(i+1, j+1)$ represents the attention weight from $(\zeta(w))_i$ to $w_j$ when generating $(\zeta(w))_k$. The $+1$ shift arises from the default \textit{bos} prefix added to every input and output sequence. Darker colors correspond to stronger attention.}
    \label{fig:xatt_example}
\end{figure}

While there are various clues that the \textit{Minimal Dyck Transformer} is implementing an interpretable algorithm, the most significant ones are the attention matrices throughout the generation of output tokens. 
An example is shown in \Cref{fig:xatt_example}.

As mentioned in \Cref{s:models}, we can obtain an explicit description of the importance of each $w_i$ to each $\zeta(w)_j$ by examining cross-attention. 
As the transformer model generates the zeta map outputs step by step, one can dynamically track which input positions receive attention at each generation step.

Across all examples, three prevalent attention patterns emerge. 
First, most prominently observed in the Step 1 attention visualization, \textit{Minimal Dyck Transformer} tends to focus on $w_i$ with the highest level $\ell_i$. 
Second, there are consistent column gaps in the attention matrices, indicating that the model ignores the $i$-th input position whenever $w_i = 1$ throughout generation. 
Third, as more output steps are generated, a triangular shape appears in the attention visualizations for many examples, suggesting a gradual shift of focus toward consecutive input positions during generation.

To validate these observations, we conducted the following analyses on a trained \textit{Minimal Dyck Transformer} with 100\% accuracy\footnote{These analyses are considered standard methods in mechanistic interpretability~\cite{sharkey2025open}.}. 
\begin{enumerate}
    \item Ablation studies: we manually removed all attention to $w_i$ if $w_i = 1$ in the model. The resulting accuracy decreased by only 1.5\%, confirming that the model indeed does not attend to such positions during generation.
    \item Probing studies: we extracted encoder hidden states and verified that the \textit{Minimal Dyck Transformer} implicitly computes the level sequence $\ell$. 
\end{enumerate}

The model’s tendency to attend to $w_i$ with the highest level $\ell_i$ and to shift attention consecutively resembles the behavior of the sweep map. 
However, its consistent disregard for positions where $w_i = 1$ suggests that it is likely learning a variant of that map.

\subsection{Scaffolding map}
From the attention patterns and subsequent analysis, we derived the following procedure, which appears to be implemented implicitly by the transformer model; we refer to it as the \textit{Scaffolding Map}.
\begin{enumerate}
    \item Start with the binary sequence representation $w = (w_i)_{i=1,\ldots,2n}$ of length $2n$ for the Dyck path over the alphabet $\{0, 1\}$.
    \item Compute levels $\ell$: start with $\ell_0 = 0$, and for each $i = 1,\dots,2n$,
    \[
      \ell_{i+1} = 
      \begin{cases}
        \ell_i + 1, & \text{if } w_i = 1,\\[2pt]
        \ell_i - 1, & \text{if } w_i = 0.
      \end{cases}
    \]

    \item Collect the right-step positions:
    \[
      R = \{\, i \in \{1,\dots,2n\} : w_i = 0 \,\}.
    \]

    \item Collect \emph{peaks} by level: for every peak position $i$ with $(w_i, w_{i+1}) = (1,0)$,
    add $i$ to the list corresponding to level $\ell_i$.

    \item Initialize:
    \[
      \textit{out} = [\,], \qquad
      \textit{agents} = [\,], \qquad
      \textit{current level} = \max\{\ell_i : i \text{ is a peak}\}.
    \]

    \item While $|\textit{out}| < 2n$:
    \begin{enumerate}
        \item Let $\textit{queue}$ be the union of
        \begin{itemize}
            \item all positions $i$ currently in $\textit{agents}$, and
            \item all peak positions $j$ in \textit{peaks} with $\ell_j = \textit{current level}$.
        \end{itemize}
        
        \item Sort $\textit{queue}$ in decreasing order and, for each $i \in \textit{queue}$ in this order, append $w_i$ to $\textit{out}$.

        \item Update all agents simultaneously as follows. For each $i$ in \textit{agents}:
        \begin{itemize}
            \item If $i \in R$, replace $i$ by $i+1$ if $i+1 \in R$ and $i+1 \le 2n$; otherwise remove $i$ from \textit{agents}. 
            \item If $i \notin R$, replace $i$ by $i-1$ if $i-1 \notin R$ and $i-1 \ge 1$; otherwise remove $i$ from \textit{agents}.
        \end{itemize}
        
        \item For each peak $j$ with $\ell_j = \textit{current level}$, attempt to append to \textit{agents} the positions
        $j+1$ if $j+1 \in R$ and $j+1 \le 2n$, and $j-1$ if $j-1 \notin R$ and $j-1 \ge 1$ (avoiding duplicates if desired).
        
        \item Decrease $\textit{current level}$ by $1$.
    \end{enumerate}

    \item Return $\textit{out}$.
\end{enumerate}

We name it the Scaffolding Map because the $1$-moves (the $w_i$'s such that $w_i = 1$) act only to build up the levels $\ell_i$, whereas the full set of levels and Dyck word entries can be reconstructed from the $0$-steps, analogous to how a physical scaffolding defines the shape of a structure. 
The elements of \textit{agents} are then imagined as positions of workers moving along this scaffolding, gathering the corresponding level and Dyck word data.

Although prior work has reverse-engineered the underlying algorithms that deep learning models implement for modular arithmetic and group operations \cite{nanda2023progress, chughtai2023toy, zhong2023clockpizzastoriesmechanistic}, to the best of our knowledge, the Scaffolding Map is the first new sequential-output algorithm discovered through interpreting a neural network for a mathematical problem.
Echoing the observations of \cite{zhong2023clockpizzastoriesmechanistic}, we also find a surprising diversity of solutions, suggesting that additional algorithmic variants remain to be uncovered. 

\printbibliography

\end{document}